\documentclass[twosided,reqno]{amsart}

\usepackage{amsmath}
\usepackage{amssymb}
\usepackage{amsthm}

\theoremstyle{theorem}
\newtheorem{theorem}{\scshape Theorem }[section]
\newtheorem{lemma}[theorem]{\scshape Lemma}

\newtheorem{proposition}[theorem]{\scshape Proposition}
\theoremstyle{definition}

\newcommand{\ma}{\mathbb}
\newcommand{\be}{\begin{equation}}
\newcommand{\ee}{\end{equation}}
\newcommand{\ben}{\begin{equation*}}
\newcommand{\een}{\end{equation*}}
\newcommand{\fa}{\frac}
\newcommand{\la}{\label}
\newcommand{\Z}{\mathbb{Z}}
\newcommand{\ck}{C_{n}^{(k)}}
\newcommand{\C}{\mathbb{C}_n^{(k)}}
\newcommand{\ch}{\widehat{C}_n^{(k)}}
\newcommand{\cn}{C_{n}}
\newcommand{\nk}{C_{n,k}}
\newcommand{\F}{\mathcal{F}}
\newcommand{\nm}{{n \brack m}}

\newcommand{\st}{S_{2}(m,n)}
\newcommand{\sn}{S_n(x)}
\newcommand{\dx}{dx_1 \ldots dx_k}
\newcommand{\ty}{\infty}
\newcommand{\I}{\int_{0}^{1}}
\newcommand{\U}{\sum_{n=0}^{\infty}}
\newcommand{\bt}{\begin{theorem}}
\newcommand{\et}{\end{theorem}}
\newcommand{\bi}{\binom}
\newcommand{\bp}{\begin{proposition}}
\newcommand{\ep}{\end{proposition}}

\numberwithin{equation}{section}

\begin{document}

\title{Higher-order Cauchy numbers and polynomials}

\author{Dae San Kim}
\address{Department of Mathematics, Sogang University, Seoul 121-742, Republic of Korea.}
\email{dskim@sogang.ac.kr}

\author{Taekyun Kim}
\address{Department of Mathematics, Kwangwoon University, Seoul 139-701, Republic of Korea}
\email{tkkim@kw.ac.kr}

\maketitle

\begin{abstract}
Recently, Komatsu introduced the concept of poly-Cauchy numbers and polynomials which generalize Cauchy numbers and polynomials.  In this paper, we consider the new concept of higher-order Cauchy numbers and polynomials which  generalize Cauchy numbers and polynomials in different direction and investigate some properties of those new class of numbers and polynomials. From our investigation, we derive some identities involving higher-order Cauchy numbers and polynomials, which generalize some relations between two kinds of Cauchy polynomials and some identities for Cauchy numbers and Stirling numbers.
\end{abstract}

\section{Introduction}

In the book of Comtet\cite{1}, two kinds of Cauchy numbers are introduced: The first kind is given by
\be\la{1}
C_n=\int_0^1 (x)_n dx,~(n\in \ma{Z}_{\geq 0})
\ee
and the second kind is given by
\be\la{2}
\widehat{C}_n=\int_0^1 (-x)_n dx, ~(n\in \ma{Z}_{\geq 0}),
\ee
where $(x)_n=x(x-1)\ldots (x-n+1)$.\\

In [2,6,7], Komatsu introduced two kinds of poly-Cauchy numbers: The poly-Cauchy numbers of the first kind $\C$ as a generalization of the Cauchy numbers are given by
\be\la{3}
\C=\I \cdots \I(x_1x_2\cdots x_k)_n dx_1dx_2\cdots dx_k,
\ee
and the poly Cauchy numbers of the second kind $\widehat{\ma{C}}_n^{(k)}$ are given by

\be\la{4}
\widehat{\ma{C}}_n^{(k)}=\I \cdots \I(-x_1x_2\cdots x_k)_n dx_1\cdots dx_k,~(n \in \ma{Z}_{\geq 0}, k \in \ma{N}).
\ee

The (signed) Stirling number of the first kind is defined by

\be\la{5}
(x)_n=\sum_{l=0}^{n}S_1(n,l)x^l,~(n \in \Z_{\geq 0}).
\ee

From (\ref{5}), we have

\be\la{6}
(\log(1+t))^n=n!\sum_{l=n}^{\ty}S_1(l,n)\fa{t^l}{l!}.
\ee

The Stirling number of the second kind is defined by the generating function to be

\be\la{7}
(e^t-1)^n=n!\sum_{l=n}^{\ty}S_2(l,n)\fa{t^l}{l!},~(see[3,4,9]).
\ee

From (\ref{1}) and (\ref{5}), we note that

\be\la{8}
\cn=\sum_{m=0}^{n}S_1(n,m)\fa{1}{m+1}=(-1)^n \sum_{m=0}^{n} {n \brack m} \fa{(-1)^m}{m+1},~(see~[1,10]),
\ee
where ${n \brack m}$ are the (unsigned) Stirling number of the first kind, arising as coefficients of the rising factorial
\ben
x^{(n)}=x(x+1)\cdots(x+n-1)=\sum_{m=0}^{n}\nm x^m,~(see~[1,8,10]).
\een

An explicit formula for $\C$ is given by
\ben
\C=(-1)^n \sum_{m=0}^{n}\nm\fa{(-1)^m}{(m+1)^k},~(n\geq0,~k\geq1),
\een

and

\ben
\widehat{\ma{C}}_n^{(k)}=(-1)^n\sum_{m=0}^{n}\nm\fa{1}{(m+1)^k},~(n\geq0,~k\geq1),~(see~[4,5,6]).
\een

The poly-Cauchy polynomials of the first kind $\C(z)$ are defined by

\be\la{9}
\C(z)=\I \cdots \I(x_1x_2 \cdots x_k-z)_m dx_1\cdots dx_k
\ee
and are expressed explicitly in terms of Stirling numbers of the first kind:

\be\la{10}
\C(z)=\sum_{m=0}^{n}\nm(-1)^{n-m}\sum_{i=0}^{m}\bi{m}{i}\fa{(-z)^i}{(m-i+1)^k},~(see~[4,5,7]).
\ee

The poly-Cauchy polynomials of the second kind $\widehat{\ma{C}}_n^{(k)}(z)$ are defined by

\be\la{11}
\widehat{\ma{C}}_n^{(k)}(z)=\I\cdots \I(-x_1 \cdots x_k+z)_m dx_1\cdots dx_k,
\ee
and are expressed explicitly in terms of Stirling numbers of the second kind:
\be\la{12}
\widehat{\ma{C}}_n^{(k)}(z)=\sum_{m=0}^{n}\nm(-1)^n\sum_{i=0}^{m}\bi{m}{i}\fa{(-z)^i}{(m-i+1)^k},~(see~[2,6,7]).
\ee

For $\alpha \in \ma{N}$, as is well known, the Bernoulli polynomials of order $\alpha$ are defined by the generating function to be
\be\la{13}
\left(\fa{t}{e^t-1}\right)^{\alpha}e^{xt}=\underbrace{\left(\fa{t}{e^t-1}\right)\times \cdots \times \left(\fa{t}{e^t-1}\right)}_{\alpha-times}e^{xt}=\U B_n^{(\alpha)}(x)\fa{t^n}{n!}.
\ee

When $x=0$, $B_n^{(\alpha)}=B_n^{(\alpha)}(0)$ are the Bernoulli numbers of order $\alpha$.~(see~[3,4,9]).\\
In this paper, we consider the new concept of higher-order Cauchy numbers and polynomials which generalize Cauchy numbers and polynomials and investigate some properties of those new class of numbers and polynomials.  From our investigation, we derive some identities involving higher-order Cauchy numbers and polynomials, which generalize some relations between two kinds of Cauchy polynomials and some identities for Cauchy numbers and Stirling numbers.

Finally, we introduce some identities of higher-order Cauchy polynomials arising from umbral calculus.

\section{Higher-order Cauchy polynomials}

For $k \in \ma{N}$, let us consider the Cauchy numbers of the first kind of order $k$ as follows:
\be\la{14}
\ck=\I\cdots \I(x_1+x_2+\cdots x_k)_n dx_1 \cdots dx_k,
\ee
where $n \in \ma{Z}_{\geq 0}$ and $k \in \ma{N}$.\\
Then, from (\ref{14}), we can derive the generating function of $\ck$ as follows:
\be\la{15}
\begin{split}
\U\ck \fa{t^n}{n!} &= \I\cdots \I\U \bi{x_1+ \cdots +x_k }{n} t^n dx_1 \cdots dx_k\\
&=\I\cdots \I(1+t)^{x_1+\cdots x_k}dx_1 \cdots dx_k.
\end{split}
\ee

It is easy to show that
\be\la{16}
\left(\fa{(1+t)^x}{\log(1+t)}\right)^{'}=(1+t)^x
\ee

Thus, by (\ref{16}), we get
\be\la{17}
\I (1+t)^xdx=\fa{t}{\log(1+t)}.
\ee

From (\ref{15}) and (\ref{17}), we have
\be\la{18}
\begin{split}
\U\ck\fa{t^n}{n!}&=\I\cdots \I(1+t)^{x_1+\cdots x_k}dx_1 \cdots dx_k\\
&=\left(\fa{t}{\log(1+t)}\right)^{k}.
\end{split}
\ee

It is known that
\be\la{19}
\left(\fa{t}{\log(1+t)}\right)^{n}(1+t)^{x-1}=\sum_{k=0}^{\ty}B_{k}^{(k-n+1)}(x)\fa{t^k}{k!},~(see~[1]).
\ee
Therefore, by (\ref{18}) and (\ref{19}), we obtain the following theorem.

\bt\la{t1}
For $n \geq 0$, we have
\ben
C_{n}^{(k)}=B_{n}^{(n-k+1)}(1).
\een
\et

From (\ref{1}), we have
\be\la{20}
\U C_n \fa{t^n}{n!}=\I \U \bi{x}{n}t^n dx = \I (1+t)^xdx=\fa{t}{\log(1+t)}.
\ee

Thus, by (\ref{18}) and (\ref{20}), we get
\be\la{21}
\U \ck \fa{t^n}{n!}=\U\left( \sum_{l_1+ \cdots +l_k=n}\bi{n}{l_1, \cdots, l_k}C_{l_1}\cdots C_{l_k}\right)\fa{t^n}{n!}.
\ee
From (\ref{14}), we note that
\be\la{22}
\begin{split}
\ck &=\I\cdots \I(x_1+\cdots +x_k)_n dx_1\cdots dx_k\\
&=\sum_{l=0}^{n}S_1(n,l)\I\cdots \I(x_1+\cdots +x_k)^l dx_1\cdots dx_k\\
&=\sum_{l=0}^{n}\sum_{l_1+ \cdots +l_k=l}S_1(n,l)\bi{l}{l_1, \cdots, l_k}\I\cdots \I x_{1}^{l_1}x_{2}^{l_2}\ldots x_{k}^{l_k}dx_1\cdots dx_k\\
&=\sum_{l=0}^{n}\sum_{l_1+ \cdots +l_k=l}\bi{l}{l_1, \cdots, l_k}S_1(n,l)\fa{1}{(l_1+1)\cdots (l_k+1)}.
\end{split}
\ee
Therefore, by (\ref{21}) and (\ref{22}), we obtain the following theorem.

\bt\la{t2}
For $n \geq 0$, we have
\ben
\begin{split}
\ck&=\sum_{l_1+ \cdots +l_k=n}\bi{n}{l_1, \cdots, l_k}C_{l_1}\cdots C_{l_k}\\
&=\sum_{l=0}^{n}\sum_{l_1+ \cdots +l_k=l}\bi{l}{l_1, \cdots, l_k}S_1(n,l)\fa{1}{(l_1+1)\cdots (l_k+1)}.
\end{split}
\een
\et
From (\ref{18}), we can derive the following equations.

\be\la{23}
\begin{split}
\U\ck\fa{(e^t-1)^n}{n!}&=\fa{1}{t^k}(e^t-1)^k=\U S_2(n+k,k)\fa{k!}{(n+k)!}t^n\\
&=\U\fa{n!k!}{(n+k)!}S_2(n+k,k)\fa{t^n}{n!}\\
&=\U\fa{S_2(n+k,k)}{\bi{n+k}{n}}\fa{t^n}{n!},
\end{split}
\ee
and

\be\la{24}
\begin{split}
\U\ck\fa{1}{n!}(e^t-1)^n&=\U\ck\sum_{m=n}^{\ty}S_2(m,n)\fa{t^m}{m!}\\
&=\sum_{m=0}^{\ty}\left(\sum_{n=0}^{m}\ck S_2(m,n)\right)\fa{t^m}{m!}.
\end{split}
\ee
Therefore, by (\ref{23}) and (\ref{24}), we obtain the following theorem.

\bt\la{t3}
For $m \in \Z_{\geq 0}$, $k \in \ma{N}$, we have
\ben
\begin{split}
S_2(m+k, k)&=\bi{m+k}{m}\sum_{n=0}^{m}\ck S_2(m,n)\\
&=\bi{m+k}{m}\sum_{n=0}^{m}B_{n}^{(n-k+1)}(1)\st.
\end{split}
\een
\et

Now, we consider the higher-order Cauchy polynomials of the first kind as follows:
\be\la{25}
\ck(x)=\I \cdots \I (x_1 + \cdots +x_k-x)_n dx_1 \cdots dx_k.
\ee
Then, by (\ref{25}), we get
\be\la{26}
\begin{split}
\ck(x)&=\sum_{l=0}^{n}S_1(n,l) \I \cdots \I (x_1 +x_2 \cdots +x_k-x)^l dx_1 \cdots dx_k\\
&=\sum_{l=0}^{n}S_1(n,l)\sum_{j=0}^{l}\bi{l}{j}(-x)^{l-j}\I \cdots \I (x_1 + \cdots +x_k)^j dx_1 \cdots dx_k\\
&=\sum_{l=0}^{n}\sum_{j=0}^{l}\sum_{j_1+\cdots j_k =j}\bi{j}{j_1, \cdots, j_k}\bi{l}{j}S_1(n,l) (-x)^{l-j}\fa{1}{(j_1+1)\cdots (j_k+1)}.
\end{split}
\ee

From (\ref{25}), we can derive the generating function of $\ck(x)$ as follows:

\be\la{27}
\begin{split}
\U\ck(x)\fa{t^n}{n!}&=\I \cdots \I\U\bi{x_1 + \cdots +x_k-x}{n}t^n \dx\\
&=\I \cdots \I(1+t)^{x_1 + \cdots +x_k-x}\dx\\
&=\left(\fa{t}{\log(1+t)}\right)^k(1+t)^{-x}.
\end{split}
\ee
It is known that
\be\la{28}
\left(\fa{t}{\log(1+t)}\right)^k(1+t)^{x}=\U B_{n}^{(n-k+1)}(x+1)\fa{t^n}{n!}.
\ee

By (\ref{27}) and (\ref{28}), we get
\ben
\ck(x)=B_{n}^{(n-k+1)}(1-x).
\een
Therefore, by (\ref{26}) and (\ref{28}), we obtain the following theorem.

\bt\la{t4}
For $n \in \Z_{\geq 0}$, $k \in \ma{N}$, we have
\ben
\begin{split}
\ck(x)&=B_{n}^{(n-k+1)}(1-x)\\
&=\sum_{l=0}^{n}\sum_{j=0}^{l}\sum_{j_1+\cdots j_k =j}\bi{j}{j_1, \cdots, j_k}\bi{l}{j}S_1(n,l) \fa{(-x)^{l-j}}{(j_1+1)\cdots (j_k+1)}.
\end{split}
\een
\et
By (\ref{27}), we see that
\be\la{29}
\begin{split}
\U\ck(x)\fa{(e^t-1)^n}{n!}&=e^{-tx}\left( \fa{e^t-1}{t}\right)^k\\
&=e^{-tx}\U S_2(n+k,k)\fa{n!k!}{(n+k)!}\fa{t^n}{n!}\\
&=\left(\sum_{l=0}^{\ty}\fa{(-x)^l}{l!}t^l\right)\left(\U \fa{S_2(n+k,k)}{\bi{n+k}{n}}\fa{t^n}{n!}\right)\\
&=\sum_{m=0}^{\ty}\left\{\sum_{n=0}^{m}\fa{\bi{m}{n}}{\bi{n+k}{n}}S_2(n+k,k)(-x)^{m-n}\right\}\fa{t^m}{m!},
\end{split}
\ee
and
\be\la{30}
\begin{split}
\U\ck(x)\fa{(e^t-1)^n}{n!}&=\U\ck(x)\sum_{m=n}^{\ty}\st\fa{t^m}{m!}\\
&=\sum_{m=0}^{\ty}\left\{\sum_{n=0}^{m}\ck(x)\st\right\}\fa{t^m}{m!}.
\end{split}
\ee
Therefore, by (\ref{29}) and (\ref{30}), we obtain the following theorem.

\bt\la{t5}
For $m \in \Z_{\geq 0}$, $k \in \ma{N}$, we have
\ben
\sum_{n=0}^{m}\fa{\bi{m}{n}}{\bi{n+k}{n}}S_2(n+k,k)(-x)^{m-n}=\sum_{n=0}^{m}\ck(x)\st.
\een
\et

We now define the Cauchy numbers of the second kind of order $k$ as follows:
\be\la{31}
\widehat{C}_n^{(k)}=\I \cdots \I(-(x_1 + \cdots +x_k))_n\dx.
\ee
From (\ref{31}), we can derive the generating function of $\widehat{C}_n^{(k)}$ as follows:

\be\la{32}
\begin{split}
\U \widehat{C}_n^{(k)}\fa{t^n}{n!}&=\I \cdots \I\U \bi{-x_1 - \cdots -x_k}{n}t^n\dx\\
&=\I \cdots \I(1+t)^{-x_1 - \cdots -x_k}\dx\\
&=\left(\fa{t}{(1+t)\log(1+t)}\right)^k.
\end{split}
\ee
Thus, by (\ref{32}), we get
\be\la{33}
\begin{split}
\sum_{m=0}^{\ty}\widehat{C}_m^{(k)}\fa{(e^t-1)^m}{m!}&=\left(\fa{e^t-1}{te^t}\right)^k\\
&=\left(\sum_{l=0}^{\ty}\fa{(-k)^l}{l!}t^l\right)\left(\sum_{m=0}^{\ty}\fa{S_2(k+m,k)k!}{(k+m)!}t^m\right)\\
&=\U\left\{\sum_{m=0}^{n}\fa{\bi{n}{m}}{\bi{k+m}{m}}(-k)^{n-m}S_2(k+m,k)\right\}\fa{t^n}{n!},
\end{split}
\ee
and

\be\la{34}
\begin{split}
\sum_{m=0}^{\ty}\widehat{C}_m^{(k)}\fa{(e^t-1)^m}{m!}&=\sum_{m=0}^{\ty}\widehat{C}_m^{(k)}\sum_{n=m}^{\ty}S_2(n,m)\fa{t^n}{n!}\\
&=\U\left(\sum_{m=0}^{n}\widehat{C}_m^{(k)}S_2(n,m)\right)\fa{t^n}{n!}.
\end{split}
\ee

Therefore, by (\ref{33}) and (\ref{34}), we obtain the following theorem.

\bt\la{t6}
For $n \geq 0$, $k \in \ma{N}$, we have
\ben
\sum_{m=0}^{n}\fa{\bi{n}{m}}{\bi{k+m}{m}}S_2(k+m,k)(-k)^{n-m}=\sum_{m=0}^{n}\widehat{C}_m^{(k)}\st.
\een
\et

We also consider the higher-order Cauchy polynomials of the second kind as follows:
\be\la{35}
\widehat{C}_n^{(k)}(x)=\I \cdots \I(x-(x_1 + \cdots +x_k))_n\dx.
\ee

By (\ref{35}), we get
\be\la{36}
\begin{split}
\widehat{C}_n^{(k)}(x)&=\sum_{l=0}^{n}S_1(n,l)\I \cdots \I(-(x_1 + \cdots +x_k)+x)^l\dx\\
&=\sum_{l=0}^{n}S_1(n,l)\sum_{i=0}^{l}\bi{l}{i}x^{l-i}(-1)^i\I \cdots \I(x_1 + \cdots +x_k)^i\dx\\
&=\sum_{l=0}^{n}\sum_{i=0}^{l}S_1(n,l)\bi{l}{i}x^{l-i}(-1)^i\sum_{j_1+\cdots j_k =i}\bi{i}{j_1, \cdots, j_k}\fa{1}{(j_1+1)\cdots (j_k+1)}\\
&=\sum_{l=0}^{n}\sum_{i=0}^{l}\sum_{j_1+\cdots j_k =i}\bi{i}{j_1, \cdots, j_k}\bi{l}{i}S_1(n,l)x^{l-i}(-1)^i\fa{1}{(j_1+1)\cdots (j_k+1)}.
\end{split}
\ee
Let us consider the generating function of the higher-order Cauchy polynomials of the second kind as follow:

\be\la{37}
\begin{split}
\U\widehat{C}_n^{(k)}(x)\fa{t^n}{n!}&=\I \cdots \I\U\bi{x-(x_1 + \cdots +x_k)}{n}t^n\dx\\
&=\I \cdots \I(1+t)^{-x_1 - \cdots -x_k+x}\dx\\
&=\left( \fa{t}{(1+t)\log(1+t)}\right)^k(1+t)^x.
\end{split}
\ee
It is not difficult to show that

\be\la{37}
\begin{split}
\U\widehat{C}_n^{(k)}(x)\fa{t^n}{n!}&=\left( \fa{t}{(1+t)\log(1+t)}\right)^k(1+t)^x\\
&=\U B_n^{(n-k+1)}(x-k+1)\fa{t^n}{n!}.
\end{split}
\ee
Therefore, by (\ref{36}) and (\ref{37}), we obtain the following theorem.

\bt\la{t7}
For $n \geq 0$, $k \in \ma{N}$, we have
\ben
\begin{split}
\widehat{C}_n^{(k)}(x)&=B_n^{(n-k+1)}(x-k+1)\\
&=\sum_{l=0}^{n}\sum_{i=0}^{l}\sum_{j_1+\cdots j_k =i}\bi{i}{j_1, \cdots, j_k}\bi{l}{i}S_1(n,l)x^{l-i}\fa{(-1)^i}{(j_1+1)\cdots (j_k+1)}.
\end{split}
\een
\et

From (\ref{37}), we note that

\be\la{39}
\begin{split}
\U\widehat{C}_n^{(k)}(x)\fa{(e^t-1)^n}{n!}&=\sum_{n=0}^{\ty}\widehat{C}_n^{(k)}(x)\sum_{m=n}^{\ty}\st\fa{t^m}{m!}\\
&=\sum_{m=0}^{\ty}\left(\sum_{n=0}^{m}\widehat{C}_n^{(k)}(x)\st\right)\fa{t^m}{m!},
\end{split}
\ee
and

\be\la{40}
\begin{split}
\U\widehat{C}_n^{(k)}(x)\fa{(e^t-1)^n}{n!}&=\left(\fa{e^t-1}{t}\right)^ke^{t(x-k)}\\
&=\left(\sum_{l=0}^{\ty}\fa{k!}{(l+k)!}S_2(l+k,k)t^l\right)\left(\U(x-k)^n\fa{t^n}{n!}\right)\\
&=\sum_{m=0}^{\ty}\left(\sum_{n=0}^{m}\fa{S_2(n+k,k)k!}{(n+k)!}(x-k)^{m-n}\fa{m!}{(m-n)!}\right)\fa{t^m}{m!}\\
&=\sum_{m=0}^{\ty}\left(\sum_{n=0}^{m}S_2(n+k,k)\fa{\bi{m}{n}}{\bi{n+k}{n}}(x-k)^{m-n}\right)\fa{t^m}{m!}
\end{split}
\ee
Therefore, by (\ref{39}) and (\ref{40}), we obtain the following theorem.

\bt\la{t8}
For $m \in \Z_{\geq 0}$, $k \in \ma{N}$, we have
\ben
\sum_{n=0}^{m}\widehat{C}_n^{(k)}(x)\st=\sum_{n=0}^{m}S_2(n+k,k)\fa{\bi{m}{n}}{\bi{n+k}{n}}(x-k)^{m-n}.
\een
\et

Now, we observe that

\be\la{41}
\begin{split}
(-1)^n\fa{\ck(x)}{n!}&=(-1)^n\I \cdots \I\bi{x_1 + \cdots +x_k-x}{n}\dx\\
&=\I \cdots \I\bi{-(x_1 + \cdots +x_k)+x+n-1}{n}\dx\\
&=\sum_{m=0}^{n}\I \cdots \I\bi{-(x_1 + \cdots +x_k)+x}{m}\bi{n-1}{n-m}\dx\\
&=\sum_{m=0}^{n}\bi{n-1}{n-m}\fa{1}{m!}m!\I \cdots \I\bi{-(x_1 + \cdots +x_k)+x}{m}\dx\\
&=\sum_{m=0}^{n}\bi{n-1}{n-m}\fa{1}{m!}\widehat{C}_m^{(k)}(x)=\sum_{m=1}^{n}\bi{n-1}{n-m}\fa{1}{m!}\widehat{C}_m^{(k)}(x).
\end{split}
\ee

Therefore, by (\ref{41}), we obtain the following theorem.

\bt\la{t9}
For $n,k \in \ma{N}$,  we have
\ben
(-1)^n\fa{\ck(x)}{n!}=\sum_{m=1}^{n}\bi{n-1}{n-m}\fa{1}{m!}\widehat{C}_m^{(k)}(x).
\een
\et

By the same method of (\ref{49}), we get
\be\la{42}
\begin{split}
(-1)^n\fa{\ch(x)}{n!}&=(-1)^n\I \cdots \I\bi{-(x_1 + \cdots +x_k)+x}{n}\dx\\
&=\I \cdots \I\bi{x_1 + \cdots +x_k-x+n-1}{n}\dx\\
&=\sum_{m=0}^{n}\bi{n-1}{n-m}\I \cdots \I\bi{x_1 + \cdots +x_k-x}{m}\dx\\
&=\sum_{m=0}^{n}\bi{n-1}{n-m}\fa{1}{m!}m!\I \cdots \I\bi{x_1 + \cdots +x_k-x}{m}\dx\\
&=\sum_{m=0}^{n}\bi{n-1}{n-m}\fa{C_m^{(k)}(x)}{m!}=\sum_{m=1}^{n}\bi{n-1}{n-m}\fa{C_m^{(k)}(x)}{m!}.
\end{split}
\ee

Therefore, by (\ref{42}), we obtain the following theorem.

\bt\la{t10}
For $n,k \in \ma{N}$,  we have
\ben
(-1)^n\fa{\ch(x)}{n!}=\sum_{m=1}^{n}\bi{n-1}{n-m}\fa{C_m^{(k)}(x)}{m!}.
\een
\et

\section{Sheffer sequences associated with higher-order Cauchy numbers and polynomials}

Let $\ma{P}$ be the algebra of polynomials in a single variable $x$ over $\ma{C}$ and let $\ma{P}^{*}$ be the vector space of all linear functionals on $\ma{P}$. The action of the linear functional $L$ on the polynomial $p(x)$ is denoted by $<L|p(x)>$, and the vector space structure on $\ma{P}^{*}$ is defined by
\begin{equation*}
\langle L+M \vert p(x)\rangle=\langle L\vert p(x)\rangle + \langle M\vert p(x)\rangle,~and~\langle cL\vert p(x)\rangle=c\langle L\vert p(x)\rangle,
\end{equation*}
where $c$ is any complex constant. Let $\F$ denote the algebra of formal power series in a single variable $t$:

\begin{equation}\label{43}
\mathcal{F}=\left\{f(t)=\sum_{k=0}^{\infty}a_{k}\frac{t^{k}}{k!}\vert a_{k}\in\mathbb{C}\right\}.
\end{equation}

The formal power series $f(t)$ defines the linear functional on $\ma{P}$ by setting
\begin{equation}\label{44}
\langle f(t) \vert x^n \rangle = a_n~~\text{for all $n \geq 0$, (see [5,9])}.
\end{equation}

By (\ref{43}) and (\ref{44}), we easily get

\begin{equation*}
\langle t^{k}\vert x^{n}\rangle =n!\delta_{n,k},~for~all~n, k\geq 0,~(see~[5,9]),
\end{equation*}
where $\delta_{n,k}$ is the Kronecker's symbol.\\

Let  $f_{L}(t)=\sum_{k\geq 0}\langle L\vert x^{k}\rangle\frac{t^{k}}{k!}$. By (\ref{45}), we get
$\langle f_{L}(t)\vert x^{n}\rangle =\langle L\vert x^{n}\rangle$. So, the map $L\longmapsto f_{L}(t)$ is a vector space isomorphism from $\mathbb{P}^{*}$ onto $\mathcal{F}$. Henceforth, $\mathcal{F}$ is thought of as both the algebra of formal power series and the space of linear functionals. We call $\mathcal{F}$ the umbral algebra. The umbral calculus is the study of umbral algebra. The order $o(f(t))$ of the non-zero power series $f(t)$  is the smallest integer $k$ for which the coefficient of $t^k$ does not vanish. (see [3,5,9]). If $o(f(t))=1$ (respectively, $o(f(t))=0$), then $f(t)$ is called a delta (repectively, an invertible) series.
For $o(f(t))=1$ and $o(g(t))=0$, there exists a unique sequence $S_n(x)$ of polynomials such that $\langle g(t)f(t)^k \vert S_n(x)\rangle =n! \delta_{n,k}$, where $n,k \geq 0$. The sequence $S_n(x)$ is called the Sheffer sequence for $(g(t),f(t))$, which is denoted by $S_n(x) \sim (g(t),f(t))$ (see [3,5,9]).  For $f(t) \in \mathcal{F}$ and $p(x) \in \mathbb{P}$, we have

\begin{equation}\label{45}
\langle e^{yt}\vert p(x)\rangle =p(y),~\langle f(t)g(t)\vert p(x)\rangle =\langle g(t)\vert f(t)p(x)\rangle=\langle f(t)\vert g(t)p(x)\rangle,
\end{equation}
and

\begin{equation}\label{46}
f(t)=\sum_{k=0}^{\infty} \langle f(t)\vert x^{k}\rangle \frac{t^k}{k!},~ p(x)=\sum_{k=0}^{\infty} \langle t^{k}\vert p(x)\rangle
\frac{x^k}{k!} ~~(\text{see [4,5,9]}.
\end{equation}
From (\ref{46}), we note that

\begin{equation}\label{47}
\langle t^{k}\vert p(x)\rangle=p^{(k)}(0),~\langle 1\vert p^{(k)}(x)\rangle=p^{(k)}(0),
\end{equation}
where $p^{(k)}(0)$ denotes the $k$-th derivative of $p(x)$ at $x=0$.  Thus, by (\ref{47}), we get

\begin{equation}\label{48}
t^{k}p(x)=p^{(k)}(x)=\fa{d^kp(x)}{dx^k},~for~all~k \geq 0,~(see[4,5,9]).
\end{equation}

Let $S_n(x) \sim (g(t),f(t))$. Then we have

\begin{equation}\label{49}
\frac{1}{g(\bar{f}(t)}e^{y \bar{f}(t)}=\sum_{k=0}^{\infty}S_k(y) \frac{t^k}{k!},~\text{for all $y \in \mathbb{C}$},
\end{equation}
where $\bar{f}(t)$ is the compositional inverse of $f(t)$ with $f(\bar{f}(t))=\bar{f}(f(t))=t$.\\

For $S_n(x) \sim (g(t),f(t))$, $q_n(x) \sim (h(t),l(t))$, let
\be\la{50}
\sn=\sum_{k=0}^{n}C_{n,k}q_k(x),
\ee
then we have

\be\la{51}
\nk=\fa{1}{k!}\langle \fa{h(\bar{f}(t))}{g(\bar{f}(t))}(l(\bar{f}(t)))^k | x^n \rangle, ~(see~[2,5,9]).
\ee

From (\ref{27}), (\ref{37}) and (\ref{49}), we note that

\be\la{52}
\ck(x) \sim \left(\left(\fa{t}{1-e^{-t}}\right)^k, e^{-t}-1\right),
\ee
and

\be\la{53}
\ch(x) \sim \left(\left(\fa{te^t}{e^{t}-1}\right)^k, e^{t}-1\right).
\ee

For $S_n(x) \sim (g(t),f(t))$, as is well known, we have

\be\la{54}
f(t)S_n(x)=nS_{n-1}(x),~(see~[3,5,9]).
\ee

By (\ref{52}), (\ref{53}) and (\ref{54}), we get

\be\la{55}
nC_{n-1}^{(k)}(x)=(e^{-t}-1)\ck(x)=\ck(x-1)-\ck(x),
\ee
and
\be\la{56}
n\widehat{C}_{n-1}^{(k)}(x)=(e^{t}-1)\ch(x)=\ch(x+1)-\ch(x).
\ee
Therefore, by (\ref{55}) and (\ref{56}), we obtain the following lemma.

\begin{lemma}\la{l11}
For $n \in \Z_{\geq 0}$, $k \in \ma{N}$, we have
\ben
nC_{n-1}^{(k)}(x)=\ck(x-1)-\ck(x),~n\widehat{C}_{n-1}^{(k)}(x)=\ch(x+1)-\ch(x).
\een
\end{lemma}

From (\ref{52}), we have
\be\la{57}
\left(\fa{t}{1-e^{-t}}\right)^k\ck(x) \sim (1, e^{-t}-1),~~(-1)^nx^{(n)} \sim (1,e^{-t}-1),
\ee
where $x^{(n)}=x(x+1)\cdots (x+n-1)$.\\
Thus, by (\ref{57}), we get

\be\la{58}
\left(\fa{t}{1-e^{-t}}\right)^k\ck(x)=(-1)^nx^{(n)}=\sum_{l=0}^{n}(-1)^lS_1(n,l)x^l.
\ee

From (\ref{58}), we have
\be\la{59}
\begin{split}
\ck(x)&=\left(\fa{1-e^{-t}}{t}\right)^k \sum_{l=0}^{n}(-1)^lS_1(n,l)x^l\\
&= \sum_{l=0}^{n} \sum_{m=0}^{l}\fa{k!}{(k+m)!}(l)_m S_2(k+m, k)S_1(n,l)(-1)^{k+l+m}x^{l-m}\\
&=\sum_{l=0}^{n} \sum_{m=0}^{l}\fa{\bi{l}{m}}{\bi{k+l-m}{k}}S_2(k+l-m, k)S_1(n,l)(-1)^{k-m}x^{m}.
\end{split}
\ee

By (\ref{53}), we get
\be\la{60}
\left(\fa{te^t}{e^{t}-1}\right)^k\ch(x)\sim (1,e^t-1),~(x)_n \sim (1, e^t-1).
\ee

Thus, from (\ref{60}), we have

\be\la{61}
\begin{split}
\ch(x)&=\left(\fa{e^{t}-1}{te^t}\right)^k (x)_n=\left(\fa{e^{t}-1}{te^t}\right)^k\sum_{l=0}^{n}S_1 (n,l)x^l\\
&=e^{-kt}\sum_{m=0}^{\ty}\fa{k!}{(k+m)!}S_2(k+m, k)t^m\sum_{l=0}^{n}S_1(n,l)x^l\\
&=e^{-kt}\sum_{l=0}^{n}\sum_{m=0}^{l}\fa{k!}{(k+m)!}S_2(k+m, k)S_1(n,l)(l)_mx^{l-m}\\
&=\sum_{l=0}^{n}\sum_{m=0}^{l}\fa{\bi{l}{m}}{\bi{m+k}{m}}S_2(k+m, k)S_1(n,l)e^{-kt}x^{l-m}\\
&=\sum_{l=0}^{n}\sum_{m=0}^{l}\fa{\bi{l}{m}}{\bi{m+k}{m}}S_2(k+m, k)S_1(n,l)(x-k)^{l-m}\\
&=\sum_{l=0}^{n}\sum_{m=0}^{l}\fa{\bi{l}{m}}{\bi{k+l-m}{k}}S_2(k+l-m, k)S_1(n,l)(x-k)^{m}.
\end{split}
\ee

Therefore, by (\ref{59}) and (\ref{61}), we obtain the following theorem.

\bt\la{t12}
For $n \in \Z_{\geq 0}$, $k \in \ma{N}$, we have
\ben
\ck(x)=\sum_{l=0}^{n} \sum_{m=0}^{l}\fa{\bi{l}{m}}{\bi{k+l-m}{k}}S_2(k+l-m, k)S_1(n,l)(-1)^{k-m}x^{m},
\een
and
\ben
\ch(x)=\sum_{l=0}^{n}\sum_{m=0}^{l}\fa{\bi{l}{m}}{\bi{k+l-m}{k}}S_2(k+l-m, k)S_1(n,l)(x-k)^{m}.
\een
\et

For $\ch(x) \sim \left(\left(\fa{te^t}{e^{t}-1}\right)^k, e^t-1 \right)$, $B_n^{(\alpha)}(x)\sim \left(\left(\fa{e^t-1}{t}\right)^{\alpha}, t \right)$, $(\alpha \in \ma{N})$,

let us assume that

\be\la{62}
\ch(x)=\sum_{m=0}^{n}C_{n,m}B_m^{(\alpha)}(x).
\ee

Then, by (\ref{50}), (\ref{51}) and (\ref{62}), we get

\be\la{63}
\begin{split}
C_{n,m}&=\fa{1}{m!}\langle \left(\fa{t}{(1+t)\log(1+t)}\right)^k \left(\fa{t}{\log(1+t)}\right)^{\alpha}(\log(1+t)^m \vert x^n \rangle\\
&=\fa{1}{m!}\langle \left(\fa{t}{(1+t)\log(1+t)}\right)^{k+\alpha} (1+t)^{\alpha} \vert (\log(1+t))^m x^n \rangle\\
&=\sum_{l=0}^{n-m}\fa{(n)_{l+m}}{(l+m)!}S_1(l+m,m) \langle \left(\fa{t}{(1+t)\log(1+t)}\right)^{k+\alpha} (1+t)^{\alpha} \vert  x^{n-l-m} \rangle\\
&=\sum_{l=0}^{n-m}\bi{n}{l+m}S_1(l+m,m)\widehat{C}_{n-l-m}^{(k+\alpha)}(\alpha)\\
&=\sum_{l=0}^{n-m}\bi{n}{l}S_1(n-l,m)\widehat{C}_{l}^{(k+\alpha)}(\alpha).
\end{split}
\ee

Therefore, by (\ref{62}) and (\ref{63}), we obtain the following theorem.

\bt\la{t13}
For $n \geq 0$, $k \in \ma{N}$, we have
\ben
\ch(x)=\sum_{m=0}^{n}\left\{\sum_{l=0}^{n-m}\bi{n}{l}S_1(n-l,m)\widehat{C}_{l}^{(k+\alpha)}(\alpha)\right\}B_n^{(\alpha)}(x).
\een
\et

%%%%%%%%%%%%%%%%%%%%%%%%%%%%%%%%%%%%%%%%%%%%%%%%%%%%%%%%%%%%%%%%%%%%%%%%%%%%%%%%%%%%%%%%%%%%

\bigskip
ACKNOWLEDGEMENTS. This work was supported by the National Research Foundation of Korea(NRF) grant funded by the Korea government(MOE)\\
(No.2012R1A1A2003786 ).
\bigskip

%%%%%%%%%%%%%%%%%%%%%%%%%%%%%%%%%%%%%%%%%%%%%%%%%%%%%%%%%%%%%%%%%%%%%%%%%%%%%%%%%%%%%%%%%%%%%%%%%%%%%%%%%%%%%%%%%%%%%%%%%%%%%%%%%%%%%%%%%%%%%%%%%%%%%%\begin{thebibliography}{99}

\end{document}